\documentclass[12pt]{amsart}
\usepackage{graphicx}
\usepackage{amssymb}

\DeclareMathOperator{\Des}{Des}

\DeclareMathOperator{\Pe}{Pk}
\DeclareMathOperator{\Sh}{Sh}

\DeclareMathOperator{\Q}{\mathcal{Q}\textit{sym}}
\DeclareMathOperator{\BQ}{\mathcal{BQ}\textit{sym}}

\title[Type B quasisymmetric functions]{The Hopf algebras of type B quasisymmetric functions and peak functions}

\author[S. K. Hsiao]{Samuel K. Hsiao}

\author[T. K. Petersen]{T. Kyle Petersen}

\begin{document}
\begin{abstract}
We show that with the appropriate choice of coproduct, the type B quasisymmetric functions form a Hopf algebra, and the recently introduced type B peak functions form a Hopf subalgebra.
\end{abstract}

\maketitle

In this note we show that the type B quasisymmetric functions form a Hopf algebra in a natural way, and that the type B peak functions of \cite{Petersen2} are a Hopf subalgebra. A future paper will explore more properties of these algebras and make connections with other results. Here we simply establish existence. The reader familiar with quasisymmetric functions and Stembridge's peak functions can skip to section \ref{sec:Bqsym}.

\section{Quasisymmetric functions}

The ring of quasisymmetric functions is well-known (see \cite{Stanley}, ch. 7.19). Recall that a quasisymmetric function is a formal series \[Q(x_1, x_2, \ldots ) \in \mathbb{Z}[[x_1, x_2,\ldots ]] \] of bounded degree such that for any composition $\alpha = (\alpha_1, \alpha_2, \ldots, \alpha_k)$, the coefficient of $x_{1}^{\alpha_1} x_{2}^{\alpha_2} \cdots x_{k}^{\alpha_k}$ is the same as the coefficient of $x_{i_1}^{\alpha_1} x_{i_2}^{\alpha_2} \cdots x_{i_k}^{\alpha_k}$ for all $i_1 < i_2 < \cdots < i_k$. Recall that a composition of $n$, written $\alpha \models n$, is an ordered tuple of positive integers $\alpha = (\alpha_1, \alpha_2, \ldots, \alpha_k)$ such that $|\alpha| = \alpha_1 + \alpha_2 + \cdots + \alpha_k = n$. In this case we say that $\alpha$ has $k$ parts, or $l(\alpha) = k$. We can put a partial order on the set of all compositions of $n$ by refinement. The covering relations are of the form \[ (\alpha_1, \ldots, \alpha_i + \alpha_{i+1}, \ldots, \alpha_k )
  \prec (\alpha_1, \ldots, \alpha_i, \alpha_{i+1}, \ldots, \alpha_k).\] Let $\Q_n$ denote the set of all quasisymmetric functions homogeneous of degree $n$. Then $\Q := \bigoplus_{n \geq 0} \Q_n$ denotes the graded ring of all quasisymmetric functions, where $\Q_0 = \mathbb{Z}$.

The most obvious basis for $\Q_n$ is the set of \emph{monomial} quasisymmetric functions, defined for any composition $\alpha = (\alpha_1, \alpha_2, \ldots, \alpha_k) \models n$,
\[ M_{\alpha} := \sum_{i_1 < i_2 < \cdots < i_k} x_{i_1}^{\alpha_1} x_{i_2}^{\alpha_2} \cdots x_{i_k}^{\alpha_k}.\]
There are $2^{n-1}$ compositions of $n$, and hence, the graded component $\Q_n$ has dimension $2^{n-1}$ as a vector space. We can form another natural basis with the \emph{fundamental} quasisymmetric functions, also indexed by compositions,
\[ F_{\alpha} := \sum_{ \alpha \leq \beta } M_{\beta},\] since, by inclusion-exclusion we can express the $M_{\alpha}$ in terms of the $F_{\alpha}$:
\[ M_{\alpha} = \sum_{ \alpha \leq \beta } (-1)^{l(\beta) - l(\alpha)} F_{\beta}.\]
As an example, \[ F_{21} = M_{21} + M_{111} = \sum_{i < j} x_i^2 x_j + \sum_{i < j < k} x_i x_j x_k  = \sum_{i \leq j < k} x_i x_j x_k .\] When $n = 0$ we have $M_{\emptyset} = F_{\emptyset} = 1$ (the unique composition of 0 is $\emptyset$).

One interesting application for the quasisymmetric functions is as a way of studying descent classes of permutations. Compositions encode descent sets in the following way. Recall that a \emph{descent} of a permutation $\pi$ of $[n] := \{ 1,2,\ldots,n\}$ is a position $i \in [n-1]$ such that $\pi_i > \pi_{i+1}$, and that an \emph{increasing run} of a permutation $\pi$ is a maximal subword of consecutive letters $\pi_{i+1} \pi_{i+2} \cdots \pi_{i+r}$ such that $\pi_{i+1} < \pi_{i+2} < \cdots < \pi_{i+r}$. By maximality, we have that if $\pi_{i+1} \pi_{i+2} \cdots \pi_{i+r}$ is an increasing run, then $i$ is a descent of $\pi$ (if $i\neq 0$), and $i+r$ is a descent of $\pi$ (if $i+r \neq n$). For any permutation $\pi$ define the \emph{descent composition}, $C(\pi)$, to be the ordered tuple listing the lengths of the increasing runs of $\pi$. If $C(\pi) = (\alpha_1, \alpha_2, \ldots, \alpha_k)$, we can recover the descent set of $\pi$:
\[ \Des(\pi) = \{ \alpha_1, \alpha_1 + \alpha_2, \ldots, \alpha_1 + \alpha_2 + \cdots + \alpha_{k-1} \}.\] Since $C(\pi)$ and $\Des(\pi)$ have the same information, we can use them interchangeably. For example the permutation $\pi = (3,4,5,2,6,1)$ has $C(\pi) = (3,2,1)$ and $\Des(\pi) = \{ 3, 5\}$.

An alternate way to define the fundamental quasisymmetric functions is as follows:
\[ F_{C(\pi)} = F_{\Des(\pi)} = \sum_{ \substack{ i_1 \leq i_2 \leq \cdots \leq i_n \\ s \in \Des(\pi) \Rightarrow i_s < i_{s+1} }} x_{i_1} x_{i_2} \cdots x_{i_n}.\] In fact, this is the definition Gessel \cite{Gessel} used when originally defining the quasisymmetric functions. With this definition, multiplication of quasisymmetric functions has a nice description in terms of the fundamental basis. If $\sigma$ is a permutation of $[n]$ and $\tau$ is a permutation of $[n+1,n+m] := \{n+1,n+2,\ldots,n+m\}$, then
\begin{equation}\label{eq:Fmult}
F_{C(\sigma)} F_{C(\tau)} = \sum_{ \pi \in \Sh(\sigma,\tau)} F_{C(\pi)},
\end{equation}
where the sum is over the set $\Sh(\sigma,\tau)$ of all \emph{shuffles} of $\sigma$ and $\tau$, i.e., those permutations of $[n+m]$ whose restriction to the letters $[n]$ is $\sigma$ and whose restriction to $[n+1,n+m]$ is $\tau$.

It has been shown that $\Q$ is a Hopf algebra \cite{Ehrenborg, MalvenutoReutenauer} with the usual product of formal power series, counit that takes functions to their constant term, and the following coproduct, \[\Delta : \Q \to \Q \otimes \Q,\] given in terms of the monomial basis:
\[ \Delta( M_{\alpha} ) = \sum_{ \beta \gamma = \alpha} M_{\beta} \otimes M_{\gamma},\] where $\beta \gamma = ( \beta_1, \ldots, \beta_{l(\beta)}, \gamma_1, \ldots, \gamma_{l(\gamma)})$ is the \emph{concatenation} of $\beta$ and $\gamma$.

The coproduct is easy to understand in terms of the following operation. For any two alphabets $X$ and $Y$ of commuting indeterminates, let $X+Y$ denote the set $X \cup Y$ totally ordered as \[ x_1 < x_2 < \cdots < y_1 < y_2 < \cdots. \] For any quasisymmetric function $Q$, the coproduct is equivalent to the map $Q(X) \mapsto Q(X+Y)$. To be precise, we can write $Q(X+Y) = \sum R(X)S(Y)$ for some quasisymmetric functions $R$ and $S$. Then we have $\Delta(Q) = \sum R \otimes S$. As an example,
\begin{align*}
M_{21}(X+Y) & = \sum_{ i < j} x_i^2 x_j + \sum_{i,j} x_i^2 y_j + \sum_{i < j} y_i^2 y_j \\
 & = M_{21}(X) + M_{2}(X)M_1(Y) + M_{21}(Y),
\end{align*}
and thus \[ \Delta(M_{21}) = M_{21} \otimes 1 + M_2 \otimes M_1 + 1 \otimes M_{21}.\]

We know that $\Q$ is a Hopf algebra because it is a graded connected bialgebra, where being connected means $\Q_0 = \mathbb{Z}$.  There is the following nice formula for its antipode \cite{Ehrenborg, MalvenutoReutenauer}:
\begin{equation}\label{eq:antipodeM}
\mbox{(Antipode) } \quad S_{Q}( M_{\alpha}) = (-1)^{l(\alpha)}\sum_{\beta \leq \alpha} M_{\overleftarrow{\beta}},
\end{equation}
where $\overleftarrow{\beta} = (\beta_k, \beta_{k-1}, \ldots, \beta_1)$; the composition formed by writing the parts of $\beta$ backwards.

\section{Peak functions}

There is a Hopf subalgebra of $\Q$ given by Stembridge's \emph{peak functions} \cite{Stembridge} (or \emph{interior peak functions} as they're called in \cite{Petersen2}), defined as follows. First, for any $\alpha \models n$, define the more general
\begin{equation}\label{eq:KA-M}
K_{\alpha} := \sum_{ \substack{ \beta \models n \\ \alpha \leq \beta^* }} 2^{l(\beta)} M_{\beta},
\end{equation}
where $\beta^*$ is the refinement of $\beta$ obtained by replacing, for $i > 1$, every part $\beta_i \geq 2$ with the two parts $(1,\beta_i -1)$. For example, if $\beta = (2,3,1,2)$, then $\beta^* = (2,1,2,1,1,1)$.

Let $\mathbf{\Pi}_n$ denote the $\mathbb{Z}$-span of the functions $K_{\alpha}$, for all $\alpha = (\alpha_1, \alpha_2, \ldots, \alpha_k) \models n$ such that $\alpha_i > 1$ unless $i = k$. In other words, all but the last part of the composition must be greater than 1. We call such compositions \emph{peak compositions}. The related $K_{\alpha}$ are the peak functions, which are shown in \cite{Stembridge} to be $F$-positive and linearly independent, forming a basis for $\mathbf{\Pi}_n$ of dimension $f_{n-1}$, the Fibonacci number defined by $f_0 = f_1 = 1$ and $f_n = f_{n-1} + f_{n-2}$ for $n \geq 2$. We define $\mathbf{\Pi} := \bigoplus_{n \geq 0} \mathbf{\Pi}_n$ to be the graded ring of peak functions.

Peak functions get their name because, just as the $F_{\alpha}$ encode descent classes of permutations, the $K_{\alpha}$, where $\alpha$ is a peak composition, encode peak classes of permutations. An \emph{interior peak} of a permutation $\pi$ of $[n]$ is any $i \in [2,n-1]$ such that $\pi_{i-1} < \pi_i > \pi_{i+1}$. Let $\Pe(\pi)$ denote the set of all interior peaks of $\pi$. Valid peak sets are precisely those subsets $I$ of $[2,n-1]$ with the property that if $i \in I$, $i-1 \notin I$. Thus we see that peak sets correspond to the peak compositions as described above. Let $\widehat{C}(\pi)$ denote the composition corresponding to $\Pe(\pi)$. Roughly speaking, it records the distances between peaks. For example, if $\pi = (3,2,7,5,4,1,8,6)$, then $\Pe(\pi) = \{ 3, 7\}$, and $\widehat{C}(\pi) = (3,4,1)$.

Now we can write the peak functions as:
\[ K_{\widehat{C}(\pi)} = K_{\Pe(\pi)} = \sum_{\substack{ S \subset [n-1] \\ \Pe(\pi) \subset S \cup (S+1) }} 2^{|S|+1} M_{S},\]
where $S+1 = \{ i : i-1 \in S\}$ and $M_S = M_{\alpha}$ where if $S = \{ s_1 < s_2 < \cdots < s_{k-1} \}$, then $\alpha = (s_1, s_2 - s_1, \ldots, n-s_{k-1})$. This definition is the original one given by Stembridge \cite{Stembridge}; it is motivated by his theory of \emph{enriched} $P$-partitions, which gives the following formula for multiplication of peak functions analogous to \eqref{eq:Fmult}:
\begin{equation}\label{eq:Kmult}
K_{\widehat{C}(\sigma)} K_{\widehat{C}(\tau)} = \sum_{ \pi \in \Sh(\sigma,\tau) } K_{\widehat{C}(\pi)},
\end{equation}
where again the sum is over all shuffles of $\sigma$ and $\tau$.

For any composition $\alpha$, let $\widehat{\alpha}$ denote the composition formed by replacing consecutive parts of 1's with their sum added to the next part to the right: \[ ( \ldots, \alpha_{i}, \underbrace{1,1,\ldots, 1}_r, \alpha_{i+r+1}, \ldots ) \mapsto (\ldots, \alpha_i, \alpha_{i+r+1}+r, \ldots ),\] where $\alpha_i, \alpha_{i+r+1} > 1$. For example, if $\alpha = (3,1,1,3,2,1,1,1)$, then $\widehat{\alpha} = (3,5,2,3)$. The map $\alpha \mapsto \widehat{\alpha}$ is thus a surjection from all compositions onto peak compositions. Stembridge then defines the map $\Theta: \Q \to \mathbf{\Pi}$ given by \[ \Theta( F_{\alpha} ) = K_{\widehat{\alpha}}.\] By considering equations \eqref{eq:Fmult} and \eqref{eq:Kmult}, we see that $\Theta$ is a graded, surjective homomorphism of rings as shown in \cite{Stembridge}. Moreover, \cite{BergeronMSW} shows that $\Theta$ is a morphism of Hopf algebras, so $\mathbf{\Pi}$ is a Hopf algebra.

\section{Type B quasisymmetric functions}\label{sec:Bqsym}

The type B quasisymmetric symmetric functions were introduced by Chow \cite{Chow}. The easiest way to describe the ring of type B quasisymmetric functions is as the multiplicative closure of the elements of $\Q$ and a new variable $x_0$. A more motivated definition is by analogy with the ordinary quasisymmetric functions. Define a \emph{pseudo-composition} of $n$, written $\alpha \Vdash n$, to be an ordered tuple of nonnegative integers $(\alpha_1, \alpha_2, \ldots, \alpha_k)$ whose sum $|\alpha| = \alpha_1 + \cdots + \alpha_k$ is $n$, where $\alpha_1 \geq 0$, $\alpha_i > 0$ for $i > 1$. In other words, given any ordinary composition $\alpha \models n$, we have two corresponding pseudo-compositions: $\alpha$ and $0\alpha = (0, \alpha_1, \ldots, \alpha_k)$. The partial order on the set of all pseudo-compositions of $n$ is again by refinement.

Now we can define a type B quasisymmetric function to be a formal series \[Q(x_0, x_1, x_2, \ldots ) \in \mathbb{Z}[[x_0, x_1, x_2,\ldots ]] \] of bounded degree such that for any pseudo-composition $\alpha = (\alpha_1, \alpha_2, \ldots, \alpha_k)$, the coefficient of $x_{0}^{\alpha_1} x_{1}^{\alpha_2} \cdots x_{k-1}^{\alpha_k}$ is the same as the coefficient of $x_{0}^{\alpha_1} x_{i_2}^{\alpha_2} \cdots x_{i_k}^{\alpha_k}$ for all $0 < i_2 < \cdots < i_k$. Let $\BQ_n$ denote the set of all quasisymmetric functions homogeneous of degree $n$. Then $\BQ := \bigoplus_{n \geq 0} \BQ_n$ denotes the ring of all type B quasisymmetric functions, where $\BQ_0 = \mathbb{Z}$.

As before we have a monomial and fundamental basis for $\BQ_n$. For any pseudo-composition $\alpha = (\alpha_1, \alpha_2, \ldots, \alpha_k) \Vdash n$, the monomial functions are
\[
 M_{\alpha} := \sum_{ i_2 < \cdots < i_k} x_{0}^{\alpha_1} x_{i_2}^{\alpha_2} \cdots x_{i_k}^{\alpha_k}.\]
There are $2^n$ pseudo-compositions of $n$, so the dimension of $\BQ_n$ is $2^n$. The fundamental basis is
\[ F_{\alpha} := \sum_{ \alpha \leq \beta } M_{\beta},\] and by inclusion-exclusion we can express the $M_{\alpha}$ in terms of the $F_{\alpha}$:
\[ M_{\alpha} = \sum_{ \alpha \leq \beta } (-1)^{l(\beta) - l(\alpha)} F_{\beta}.\]
As an example, \[ F_{21} = M_{21} + M_{021} + M_{111} + M_{0111} = \sum_{0 \leq i \leq j < k} x_i x_j x_k.\] A special case is $n = 0$, where we define $M_{\emptyset} = M_0 = F_{\emptyset} = F_0 = 1$. We further note that if $\alpha_1 = 0$, then $M_{\alpha}$ and $F_{\alpha}$ have no $x_0$'s, i.e., they are ordinary quasisymmetric functions.

Pseudo-compositions can encode descent classes of signed permutations. If $\pi$ is any signed permutation of $[n]$, the pseudo-composition $C_B(\pi) \Vdash n$ lists the lengths of the increasing runs of $\pi$, with the extra rule that the first part of $C_B(\pi)$ is zero if $\pi_1 < 0$. This corresponds to type B descent sets, since \[ \Des_B(\pi) := \{ i \in [0,n-1] : \pi_i > \pi_{i+1} \},\] where $\pi_0 = 0$. Again, we have an easy correspondence between the two. If $C_B(\pi) = (\alpha_1, \ldots, \alpha_k)$, \[ \Des_B(\pi) = \{ \alpha_1, \alpha_1 + \alpha_2, \ldots, \alpha_1 + \alpha_2 + \cdots + \alpha_{k-1} \}.\] Two examples are $\sigma = (1,3,-2,4)$ and $\tau = (-3,-2,4,1)$, for which $C_B(\sigma) = (2,2)$, $\Des_B(\sigma) = \{ 2\}$, $C_B(\tau) = (0,3,1)$, and $\Des_B(\tau) = \{ 0, 3\}$.

The fundamental basis for the type B quasisymmetric functions encode descent classes of signed permutations just as in the type A case:
\[ F_{C_B(\pi)} = F_{\Des_B(\pi)} = \sum_{ \substack{ 0 \leq i_1 \leq \cdots \leq i_n \\ s \in \Des_B(\pi) \Rightarrow i_s < i_{s+1} }} x_{i_1} \cdots x_{i_n}.\] Multiplication in $\BQ$ can then be described by the following formula:
\begin{equation}\label{eq:typeBshuf}
F_{C_B(\sigma)} F_{C_B(\tau)} = \sum_{ \pi \in \Sh(\sigma,\tau) } F_{C_B(\pi)},
\end{equation}
the sum again over all shuffles.

The main goal of this paper is to show that $\BQ$ is a Hopf algebra. We take the product on $\BQ$ to be the usual product of formal power series and the counit again maps type B quasisymmetric functions to their constant term. We will describe the coproduct on $\BQ$ shortly, but first let us recall Chow's coproduct \cite{Chow} on $\BQ$ that makes $\BQ$ a module coalgebra over $\Q$. Define $\Delta': \BQ \to \BQ \otimes \Q$ by:
\[
\Delta'(M_{\alpha}) = \sum_{ \beta \gamma = \alpha } M_{\beta} \otimes M_{0\gamma},
\]
where if $\gamma_1 = 0$, we understand $0 \gamma = \gamma$. This coproduct arises naturally if we consider $X_0 + Y$, the union of sets $X_0 := \{x_0\} \cup X$ and $Y$, with total order \[ x_0 < x_1 < x_2 < \cdots < y_1 < y_2 < \cdots, \] (note that there is no $y_0$). Then, for example,
\begin{align*}
 M_{12}(X_0 + Y) & = \sum_{0<i} x_0 x_i^2 + \sum_{0<i}x_0 y_i^2 + \sum_{0<i < j} y_i y_j^2 \\
 &= M_{12}(X_0) + M_{1}(X_0)M_{02}(Y) + M_{012}(Y).
\end{align*}
If $\Delta$ is the coproduct on $\Q$, $m$ denotes multiplication of power series, and $I$ is the identity, Chow \cite{Chow} shows that $\Delta'$ satisfies \[m(I \otimes \Delta)\Delta' = m(\Delta' \otimes I)\Delta',\] which shows that $\BQ$ is a module coalgebra over $\Q$.

To get the Hopf algebra structure we need a new coproduct, defined as follows. Let $\alpha$ be any ordinary composition. Then $\Delta: \BQ \to \BQ \otimes \BQ$ is:
\begin{equation}
\Delta(M_{n\alpha}) = \sum_{\beta\gamma = \alpha} \sum_{i=0}^n \binom{n}{i} M_{i\beta} \otimes M_{(n-i)\gamma}, \label{eq:cop}
\end{equation}
with $\Delta(M_{\emptyset}) = \Delta(M_0) = 1 \otimes 1$. For $n = 0$, note that rule \eqref{eq:cop} is the same as for ordinary quasisymmetric functions (this makes sense as $M_{0\alpha} \in \Q$). As an example where $n \neq 0$, we have \[ \Delta(M_{121}) = M_{121} \otimes 1 + M_{021} \otimes M_1 + M_{12} \otimes M_{01} + M_{02} \otimes M_{11} + M_1 \otimes M_{021} + 1 \otimes M_{121}.\] Another way of understanding the coproduct is by taking $X_0 + Y_0$, the union of sets $X_0$ and $Y_0$. We define the map $M_{\alpha}(X_0) \mapsto M_{\alpha}(X_0 + Y_0)$ as for ordinary quasisymmetric functions, with the caveat that $x_0 \mapsto x_0 + y_0$ directly. For example,
\begin{align*} M_{211}(X_0 + Y_0) & = (x_0 + y_0)^2 \Big( \sum_{0 < i < j} x_i x_j + \sum_{0< i, j} x_i y_j + \sum_{0< i < j } y_i y_j \Big) \\
 & = \sum_{ 0 < i < j } ( x_0^2 x_i x_j + 2 x_0 x_i x_j y_0 + x_i x_j y_0^2 ) \\
 & \quad + \sum_{0 < i,j} ( x_0^2 x_i y_j + 2 x_0 x_i y_0 y_j + x_i y_0^2 y_j) \\
 & \quad + \sum_{0 < i < j} (x_0^2 y_i y_j + 2 x_0 y_0 y_i y_j + y_0^2 y_i y_j) \\
 & = M_{211}(X_0) + 2 M_{111}(X_0)M_{01}(Y_0) + M_{011}(X_0)M_{2}(Y_0) \\
 & \quad + M_{21}(X_0)M_{01}(Y_0) + 2M_{11}(X_0)M_{11}(Y_0) + M_{01}(X_0)M_{21}(Y_0) \\
 & \quad + M_{2}(X_0)M_{011}(Y_0) + 2M_1(X_0)M_{111}(Y_0) + M_{211}(Y_0),
\end{align*}
corresponding to:
\begin{align*} \Delta(M_{211}) & = M_{211} \otimes 1 + 2 M_{111} \otimes M_{1} + M_{011} \otimes M_{2} \\
 & \quad + M_{21} \otimes M_{01} + 2 M_{11} \otimes M_{11} + M_{01} \otimes M_{21} \\
 & \quad + M_{2} \otimes M_{011} + 2 M_1 \otimes M_{111} + 1\otimes M_{211}.
\end{align*}

We will now show that $\BQ$ is a Hopf algebra under this coproduct. Let $\mathbb{Z}[x_0] = \mathbb{Z}[M_1]$ denote the ring of polynomials in one variable. Notice that we have a natural isomorphism of algebras \[\BQ \cong \Big( \mathbb{Z}[M_1] \otimes \Q\Big)\] given, for any ordinary composition $\alpha$ and any $i \geq 0$ by \[ M_{i\alpha} \mapsto M_i \otimes M_{0\alpha} = M_1^i \otimes M_{0\alpha}, \] with inverse mapping given by the usual product: \[ M_i \otimes M_{0\alpha} \mapsto M_i \cdot M_{0\alpha} = M_{i\alpha}.\] The Hopf structure on $\Q$ has already been described; let its coproduct be denoted $\Delta_2$. There is a Hopf structure on $\mathbb{Z}[M_1]$ given by the usual product of polynomials and coproduct $\Delta_1: \mathbb{Z}[M_1] \to \mathbb{Z}[M_1] \otimes \mathbb{Z}[M_1]$ given by \[M_1 \mapsto M_1 \otimes 1 + 1 \otimes M_1,\] and $\Delta_1(1) = 1 \otimes 1$. In other words, a polynomial $f(x_0)$ maps to the two-variable polynomial $f(x_0 + y_0)$. With this Hopf structure $\mathbb{Z}[M_1]$ is isomorphic to the binomial Hopf algebra, a well-studied object especially in the context of the umbral calculus \cite{Ray}.

It is straightforward to verify that $\mathbb{Z}[M_1] \otimes \Q$ is a Hopf algebra with coproduct $\Delta = \Delta_1 \otimes \Delta_2$. Translated in terms of $\BQ$, this coproduct is just \eqref{eq:cop} as defined above. Combining \eqref{eq:antipodeM} with the well-known formula for the antipode of the binomial Hopf algebra, $S(M_n) = (-1)^n M_n$, we get the following formula:
\[\text{(Antipode on $\BQ$) } \quad S_{B}(M_{n\alpha}) = (-1)^{l(\alpha)+n}\sum_{\beta \leq \alpha} M_{n\overleftarrow{\beta}},
\]
where $\overleftarrow{\beta}$ is the composition formed by writing the parts of $\beta$ backwards.

\section{Type B peak functions}

For any pseudo-composition $\alpha \Vdash n$, we define the following type B quasisymmetric function in $\BQ_n$,
\begin{equation}\label{eq:KB-M}
K_{\alpha} := \sum_{ \substack{ \beta \Vdash n \\ \alpha \leq \beta^* }} 2^{l(\beta)-1} M_{\beta},
\end{equation}
 where as before $\beta^*$ is the refinement of $\beta$ obtained by replacing, for $i > 1$, every part $\beta_i \geq 2$ with the two parts $(1,\beta_i -1)$. Notice that if $\alpha_1 = 0$, then $K_{\alpha}$ is a type A peak function.

Let $\mathbf{\Pi}_{B,n}$ denote the $\mathbb{Z}$-span of the functions $K_{\alpha}$, for all $\alpha = (\alpha_1, \alpha_2, \ldots, \alpha_k) \Vdash n$ such that the ``middle parts" $\alpha_2, \ldots, \alpha_{k-1}$ are all greater than 1. We call such compositions \emph{type B peak pseudo-compositions}. Define the ring $\mathbf{\Pi}_B := \bigoplus_{n \geq 0} \mathbf{\Pi}_{B,n}$, which we call the \emph{type B peak functions}. In \cite{Petersen2} it is shown that for type B peak pseudo-compositions $\alpha$, the $K_{\alpha}$ are $F$-positive and linearly independent. The graded components $\mathbf{\Pi}_{B,n}$ has dimension $f_{n+1}$, where $f_n$ is the $n$-th Fibonacci number as defined previously.

The type B peak functions are first introduced and discussed in \cite{Petersen2}. The \emph{type B peak set} of a signed permutation $\pi$ of $[n]$ is defined as \[ \Pe_B(\pi) := \{ i \in [0,n-1] : \pi_{i-1} < \pi_i > \pi_{i+1} \},\] where $\pi_0 = 0, \pi_{-1} = -\infty$. In other words, the type B peak set is the ordinary peak set, along with 1 if $0 < \pi_1 > \pi_2$ and along with $0$ if $0 > \pi_1$. Let $\widehat{C}_B(\pi)$ be the pseudo-composition corresponding to $\Pe_B(\pi)$. As examples, if $\sigma = (-3,2,-4,5,1)$, $\tau = (3,-2,-1,5,4)$, then $\Pe_B(\sigma) = \{0,2,4\}$ and $\widehat{C}_B(\sigma) = (0,2,2,1)$, $\Pe_B(\tau) = \{1,4\}$ and $\widehat{C}_B(\tau) = (1,3,1)$.

Expressions similar to \eqref{eq:KA-M} and \eqref{eq:Kmult} hold for signed permutations:
\begin{equation}\label{eq:BKdef} K_{\widehat{C}_B(\pi)} = K_{\Pe_B(\pi)} = \sum_{ \substack{ S \subset [0,n-1] \\ \Pe_B(\pi) \subset S \cup (S+1) }} 2^{|S|} M_{S},
\end{equation}
\begin{equation}\label{eq:Bpeakshuf}
K_{\widehat{C}_B(\sigma)} K_{\widehat{C}_B(\tau)} = \sum_{ \pi \in \Sh(\sigma,\tau)} K_{\widehat{C}_B(\pi)},
\end{equation}
where $\pi$, $\sigma$, and $\tau$ are now signed permutations.

For any pseudo-composition $\alpha = (\alpha_1, \alpha_2, \ldots, \alpha_k)$, define $\widehat{\alpha} = \alpha_1 \widehat{(\alpha_2, \ldots, \alpha_k)}$, where we recall that $\widehat{\beta}$ is the peak composition formed from an ordinary composition $\beta$ by replacing consecutive parts of size 1 with their sum, added to the part to their right. Then the map $\alpha \mapsto \widehat{\alpha}$ is a surjection from all pseudo-compositions onto peak pseudo-compositions. For example, if $\alpha = (1,1,3,2,1,1,3,1)$, then $\widehat{\alpha} = (1,4,2,5,1)$.

We can then define the map $\Theta_B: \BQ \to \mathbf{\Pi}_{B}$ by
\begin{equation}\label{eq:ThetaB-F}
 \Theta_B (F_{\alpha}) = K_{\widehat{\alpha}},
\end{equation}
which by equations \eqref{eq:typeBshuf} and \eqref{eq:Bpeakshuf} defines a graded surjective homomorphism of rings. Notice that the restriction of $\Theta_B$ to $\Q$ agrees with Stembridge's map $\Theta$.

We now want to show that $\Theta_B: \BQ \to \BQ$ is a homomorphism of Hopf algebras since it will follow that $\mathbf{\Pi}_B$, as the image of $\Theta_B$, is a Hopf subalgebra of $\BQ$. A quick calculation yields  \[\Theta_B(M_1) = M_1.\] Therefore $\Theta_B$ corresponds to the map $e \otimes \Theta$ under the natural isomorphism $\BQ \cong \mathbb{Z}[M_1] \otimes \Q$, where $e$ denotes the identity map on $\mathbb{Z}[M_1]$. Since $e$ and $\Theta$ are both Hopf homomorphisms, so is $\Theta_B$.

\section{Connections}

The algebra $\BQ$ is a very natural object, in the sense that it is isomorphic to the tensor of two fundamental algebras, namely, the binomial algebra and the quasisymmetric functions. In future work, we hope to study the algebras $\BQ$ and $\mathbf{\Pi}_B$ in more detail. We would like to derive another coproduct formula for $\BQ$ in terms of the fundamental basis and a coproduct formula for the type B peak algebra (in terms of any basis). Further, it would be nice to make connections with Hopf algebras of posets and with the theory of combinatorial Hopf algebras \cite{AguiarBergeronSottile}.

\end{document}